\documentclass[12pt,reqno, oneside]{amsart}

\newtheorem{theorem}{Theorem}

\newcommand{\egf}[1] {\sum_{n=0}^\infty #1\frac{x^n}{n!}}
\newcommand{\sie}{symmetric inclusion-exclusion}
\newcommand{\aie}{asymmetric inclusion-exclusion}
\newcommand{\E}{\overline{E}}

\renewcommand{\S}{\overline{S}}
\newcommand{\T}{\overline{T}}

\renewcommand{\P}{\mathbf P}
\newcommand{\Above}[2]{\genfrac{}{}{0pt}{}{#1}{#2}}

\let \axx=\alpha
\let\bxx=\beta
\let\rx=r

\begin{document}

\newbox\Adr
\setbox\Adr\vbox{
\centerline{\sc Ira M. Gessel$^*$}
\vskip7pt
\centerline{Department of Mathematics, MS 050}
\centerline{Brandeis University}
\centerline{Waltham, MA 02454-9110}
\centerline{\tt gessel@brandeis.edu}
}

\title[Symmetric inclusion-exclusion]{Symmetric inclusion-exclusion}

\author[Ira M. Gessel]{\box\Adr} 
\address{Department of Mathematics\\
    Brandeis University\\
    Waltham, MA 02454-9110}
 \email{gessel@brandeis.edu}
\date{July 23, 2005}
\thanks{$^*$Partially supported by NSF Grant DMS-0200596}
\begin{abstract}
 One form of the inclusion-exclusion principle asserts that if $A$ and $B$ are functions of finite sets then the formulas
$A(S) = \sum_{T\subseteq S}B(T)$ and
$B(S) = \sum_{T\subseteq S}(-1)^{|S|-|T|}A(T)$ are equivalent. If we replace $B(S)$ by $(-1)^{|S|}B(S)$ then these formulas take on the symmetric form
\begin{align*}
A(S) &= \sum_{T\subseteq S}(-1)^{|T|}B(T)\\
B(S) &= \sum_{T\subseteq S}(-1)^{|T|}A(T).
\end{align*}
which we call \emph{symmetric inclusion-exclusion}. We study instances of symmetric inclusion-exclusion in which the functions $A$ and $B$ have combinatorial or probabilistic interpretations. In particular, we study cases related to the P\'olya-Eggenberger urn model in which  $A(S)$ and $B(S)$ depend only on the cardinality of $S$.
\end{abstract}
 
 \vglue -30pt
 \maketitle

\section{Inclusion-exclusion}

Let $A$ and $B$ be two functions defined on 
a set $D$ of finite sets.  We assume that if $S\in D$ and  $T\subseteq S$ then $T\in D$. Then the inclusion-exclusion principle asserts that  the following are equivalent:
\begin{equation}
\label{e-1}
\begin{aligned}
A(S) &= \sum_{T\subseteq S}B(T).\\
B(S) &= \sum_{T\subseteq S}(-1)^{|S|-|T|}A(T).
\end{aligned}
\end{equation}

The special case in which $A(S)$ and $B(S)$ depend only on $|S|$ is especially important: for two sequences $(A_n)$ and $(B_n)$, the following are equivalent:
\begin{equation}
\label{e-2}
\begin{aligned}
A_n &= \sum_{k=0}^n \binom nk B_k\\
B_n &= \sum_{k=0}^n (-1)^{n-k} \binom nk A_k
\end{aligned}
\end{equation}

It is very easy to find pairs of functions with combinatorial interpretations that satisfy \eqref{e-1}: we may choose $B$ to be an arbitrary function with a combinatorial interpretation. Then since \eqref{e-1} expresses $A$ in terms of $B$ with positive coefficients, $A$ will also have a combinatorial interpretation. Similarly, it is easy to find instances of \eqref{e-2} with combinatorial interpretations.

In terms of exponential generating functions, \eqref{e-2} gives the familiar formulas 
\begin{equation*}
\label{e-egf}
\begin{aligned}
F_A(x) &=e^xF_B(x)\\
F_B(x)&=e^{-x}F_A(x),
\end{aligned}
\end{equation*}
where
\begin{equation*}
F_A(x) = \egf{A_n} \hbox{\quad and \quad } F_B(x) = \egf{B_n}.
\end{equation*}

Now in \eqref{e-1}, let us replace $B(T)$ with $(-1)^{|T|}\bxx(T)$ 
and $A(S)$ with $\axx(S)$, and multiply the second equation by $(-1)^{|S|}$.
We obtain
\begin{equation}
\label{e-1*}
\begin{aligned}
\axx(S) &= \sum_{T\subseteq S}(-1)^{|T|}\bxx(T)\\
\bxx(S) &= \sum_{T\subseteq S}(-1)^{|T|}\axx(T).
\end{aligned}
\end{equation}
The corresponding substitution in 
\eqref{e-2} gives
\begin{equation}
\label{e-2*}
\begin{aligned}
\axx_n &= \sum_{k=0}^n (-1)^k \binom nk \bxx_k\\
\bxx_n &= \sum_{k=0}^n (-1)^k \binom nk \axx_k.
\end{aligned}
\end{equation}
In terms of exponential generating functions, \eqref{e-2*} is the unfamiliar-looking
\begin{equation}
\label{e-egf*}
\begin{aligned}
F_\axx(x) &= e^x F_\bxx(-x)\\
F_\bxx(x) &= e^x F_\axx(-x)
\end{aligned}
\end{equation}
 We will call functions or sequences $A$ and $B$ that satisfy  \eqref{e-1} and \eqref{e-2} \emph{asymmetric inclusion-exclusion pairs} and we will call instances of \eqref{e-1*} and \eqref{e-2*} \emph{symmetric inclusion-exclusion pairs}.
The  goal of this paper is to study symmetric inclusion-exclusion pairs with combinatorial or probabilistic interpretations. We shall see that for the case in which the domain of the functions $\axx$ and $\bxx$ in \eqref{e-1*}  is the set of subsets of a finite set, or if the sequences $(\axx_n)$ and $(\bxx_n)$ in \eqref{e-2*} are finite (i.e., $\axx_n$ and $\bxx_n$ are defined only for $n\le N$, for some $N$), there is a very simple construction of symmetric inclusion-exclusion pairs analogous to starting with an arbitrary function $B$ in \eqref{e-1} or with an arbitrary sequence $(B_n)$ in \eqref{e-2}. We also give a very simple probabilistic setting for~\eqref{e-1*} in which the cardinalities of the sets in the domain of $\axx$ and $\bxx$ are unbounded. Finding a probabilistic interpretation of \eqref{e-2*} with infinite sequences does not seem so straightforward. We discuss some examples, related to the P\'olya-Eggenberger urn model, of infinite sequences satisfying 
\eqref{e-2*}, but we have not found a general theory of such sequences.
 
\section{Difference Tables}
 There is a close connection between  \aie\ pairs and \sie\ pairs. We illustrate with the special case of the derangement numbers. We take $A_n=n!$, so
$B_n = \sum_{k=0}^n(-1)^{n-k}\binom nk k!$ is the $n$th derangement number.
We can compute $B_n$ from $A_n$ by using a \emph{difference table}, in which the numbers $A_n=n!$ appear in the zeroth row and each number in a row below the zeroth row is the number above it to the right minus the number above it to the left.
\bigskip
\[\vbox{\halign{
&\vrule height 0pt width 0pt depth 14pt\hbox to 15pt{\hss$#$\hss}\cr
1 && 1 && 2 && 6 && 24 && 120\cr
&0 && 1 && 4 && 18 && 96\cr
&&1 && 3 && 14 && 78\cr
&&&2 && 11 && 64\cr
&&&&9 && 53\cr
&&&&&44\cr
}}\]
Then the numbers $B_n=1,0,1,2,9,44,\dots$  appear in the zeroth diagonal. Now let us rotate this triangular array $60^\circ$ counterclockwise, obtaining

\bigskip
\[\vbox{\halign{
&\vrule height 0pt width 0pt depth 14pt\hbox to 15pt{\hss$#$\hss}\cr
&&&&&120\cr
&&&&24 && 96\cr
&&& 6 && 18 && 78\cr
&&2 && 4 &&14&&  64\cr
&1&&1&&3&&11 && 53\cr
1&&0&&1&&2&&9&&44\cr
}}\]

In this array, every number above the bottom row is the sum of the two numbers below it. Let us define the numbers $\axx_i=(5-i)!$ to be the numbers in the left diagonal and let us define the numbers $\bxx_i=120, 96, 78,\dots$ to be the numbers in the right diagonal. It is not hard to check that the sequences $(\axx_n)$ and $(\bxx_n)$ satisfy \eqref{e-2*} for $n\le 5$. We also note that these sequences have combinatorial interpretations:
$\axx_i$ is the number of permutations of $\{1,2,3,4,5\}$, in which 1, 2, \dots, $i$ (or any $i$ numbers from $\{1,2,3,4,5\}$) are all fixed points and $\bxx_i$ is the number of permutations of $\{1,2,3,4,5\}$ in which 1, 2, \dots, $i$ (or any $i$ numbers from $\{1,2,3,4,5\}$) are all nonfixed points.

The general situation is described by the following theorem,  which can be proved directly by a straightforward computation that we omit, or as a consequence of Theorem \ref{t-2} below. 

\begin{theorem}
\label{t-1}
Let $B_0, B_1, \dots, B_N$ be arbitrary and let $\axx_n$ and $\bxx_n$ for $0\le n\le N$ be defined by 
\begin{align}
\axx_{N-n}&=\sum_{k=0}^n \binom nk B_k\\
\bxx_{N-n}&=\sum_{k=0}^n \binom nk B_{N-k}.
\end{align}
Then the sequences $(\axx_n)$ and $(\bxx_n)$ satisfy \eqref{e-2*} for $0\le n\le N$.
\end{theorem}

We have a similar result that allows us to construct instances of \eqref{e-2} with combinatorial interpretations. The proof is also a straightforward verification, which we omit.

\begin{theorem}
\label{t-2}
Let $\Delta$ be a finite set, and let $B$ be a function defined on the subsets of~$\Delta$. Let the functions $\axx$ and $\bxx$ be defined on subsets of $\Delta$ by 
\begin{align*}
\axx(S)&=\sum_{T\subseteq\S} B(T)\\
\bxx(S)&=\sum_{T\subseteq\S} B(\T),
\end{align*}
where the complements are with respect to $\Delta$. Then the functions $\axx$ and $\bxx$ satisfy \eqref{e-1*}.
\end{theorem}


We note that Theorem \ref{t-1} is equivalent to the case of Theorem \ref{t-2} in which $B(T)$ depends only on the cardinality of $T$.

Theorems \ref{t-1} and \ref{t-2} allow us to construct instances of  \eqref{e-1*} and \eqref{e-2*} with combinatorial interpretations.
However, in any example constructed in this way,  the cardinalities of the sets involved are bounded and cannot in general be extended to unbounded cardinalities. 

For example, if we want to extend the triangular array given before Theorem \ref{t-1} by adding a row at the bottom consisting of nonnegative real numbers, the 0 in the bottom row forces both numbers below it be 0, and the next row is then forced to have three consecutive 0's. We then have the array
\bigskip
\[\vbox{\halign{
&\vrule height 0pt width 0pt depth 14pt\hbox to 15pt{\hss$#$\hss}\cr
&&&&&&&120\cr
&&&&&&24 && 96\cr
&&&&& 6 && 18 && 78\cr
&&&&2 && 4 &&14&&  64\cr
&&&1&&1&&3&&11 && 53\cr
&&1&&0&&1&&2&&9&&44\cr
&1&&0&&0&&1&&1&&8&&36\cr
1&&0&&0&&0&&1&&0&&8&&28\cr
}}\]
to which it is impossible to add another nonnegative row at the bottom.

In the next section we discuss a setting for \eqref{e-1*} in
which the sets may have unbounded cardinalities and the functions
$\axx$ and $\bxx$ are \emph{probabilities} rather than integers.

\section{General \sie}
\label{s-3}
We consider a probability space containing sets $E_i$, where $i$
ranges over an index set $\Delta$, which in our examples will be
the set of positive integers. Thus we have a probability function
$P$, defined on all sets generated from the $E_i$ by complements
and finite unions and intersections, with the property that if $S$
and $T$ are disjoint then $P(S\cup T) = P(S) +P(T)$. (The usual
definition of a probability function requires nonnegativity and
countable additivity, and that the probability of
the whole space is 1, but these conditions are not necessary for
our results.) For every finite subset $S\subseteq \Delta$ we set
\begin{align}
\label{e-ab}
\axx(S) &= P\Bigl(\,\bigcap_{i\in S} E_i\Bigr)\\
\bxx(S) &= P\Bigl(\,\bigcap_{i\in S} \E_i\Bigr).
\end{align}

\begin{theorem}
\label{t-3}
The functions $\axx$ and $\bxx$ defined by \eqref{e-ab} form a symmetric inclusion-exclusion pair; i.e., 
\begin{equation*}
\begin{aligned}
\axx(S) &= \sum_{T\subseteq S}(-1)^{|T|}\bxx(T)\\
\bxx(S) &= \sum_{T\subseteq S}(-1)^{|T|}\axx(T).
\end{aligned}
\end{equation*}
\end{theorem}

 \begin{proof}
By symmetry, it is enough to prove the second formula. 
Without loss of generality, we may assume that $S = \{1,2,\dots, m\}$.
A well-known form of the inclusion-principle (see, e.g., \cite[p.~6]{ross}) asserts that 
\begin{multline*}
P(E_1\cup E_2\cup \cdots \cup E_m) \\
=\sum_{1\le i \le m} P(E_i) 
\,-\sum_{1\le i<j\le m}P(E_i\cap E_j) +\cdots + (-1)^{m-1} P(E_1\cap\cdots \cap E_m).
\end{multline*}
This is equivalent to the desired formula. \end{proof}





We note that if $\Delta$ is finite, Theorem \ref{t-3} may be derived from Theorem \ref{t-2}
by taking
\begin{equation*}
B(T) = P\biggl(\Bigl(\,\bigcap_{i\in \Delta -T} E_i\Bigr)
  \bigcap\,   \Bigl(\,\bigcap_{i\in T} \E_i\Bigr)\biggr).
\end{equation*}
In the general case, we can deduce Theorem \ref{t-3} from Theorem \ref{t-2} by restricting
to a finite subset of $\Delta$ that contains $S$.

The most interesting examples of Theorem \ref{t-3} are those in
which $\axx(S)$ and $\bxx(S)$ depend only on the cardinalities of $S$.
Before discussing these examples, we give an example which is not of this
type. We consider the set of infinite sequences $(r_1,
r_2,\dots)$ of real numbers between 0 and 1, and we assign
probabilities in the obvious way to subsets defined by a finite
number of inequalities among the $r_i$. Let us take $E_i$ to be
the set of sequences satisfying $r_i<r_{i+1}$. Then for any set
$S$ of positive integers, $\axx(S)$ is the probability that
$r_i<r_{i+1}$ for all $i\in S$.  To give a formula for $\axx(S)$, we express
$S$ as a disjoint union of blocks of consecutive integers with
gaps in between: $S = V_1 \cup V_2 \cup \dots \cup V_m$, where
$V_i= \{u_i, u_i+ 1, \dots, u_i+ s_i -1\}$ and $u_{i+1}>u_i+s_i$.
Then
\[\axx(S) = \frac1{(s_1 +1)!\,(s_2+1)!\cdots (s_m+1)!},\]
and by symmetry, $\bxx(S)=\axx(S)$.
\section{A simple example}
We now give a very simple example of Theorem \ref{t-3} in which
$\axx(S)$ depends only on the cardinality of $S$, so that we have an
instance of \eqref{e-2*}. We flip a coin infinitely many times.
Each flip comes up heads with probability $p$ and tails with
probability $1-p$. For $i\in \P = \{1,2,3, \dots\}$  we let $E_i$
be the event that the $i$th flip is a head. Then for any finite
subset $S\subseteq \P$
\begin{align*}
\axx(S) &= P\Bigl(\,\bigcap_{i\in S} E_i\Bigr)=p^{|S|}\\
\bxx(S) &= P\Bigl(\,\bigcap_{i\in S} \E_i\Bigr)=(1-p)^{|S|}
\end{align*}
Then Theorem \ref{t-3} gives
\[(1-p)^{n}=\sum_{k=0}^n (-1)^k \binom nk p^k
\]
and \eqref{e-egf*} is the trivial formula
$e^{(1-p)x}=e^x e^{-px}$. 
The example of the previous section, when all of the $s_i$ are equal to 1, is equivalent to the case $p=1/2$ of this example.
\section{The P\'olya-Eggenberger urn model}
In the previous example, the events $E_i$ were independent, which makes the situation easy to analyze, but not very interesting. Here we give a more interesting case of symmetric inclusion-exclusion sequences in which the $E_i$ are not independent. It is based on a classical probabilistic model called the \emph{P\'olya} or \emph{P\'olya-Eggenberger} urn model \cite[p.~120--121]{feller}.

We have an urn that initially contains $r$ red balls and $b$ black balls. At each step we choose a ball at random from the urn. (So the probability of picking a red ball at the first step is $r/(r+b)$ and the probability of choosing a black ball is $b/(r+b)$.) We then replace the ball we have picked and add another ball of the same color to the same urn. We repeat this procedure forever.

The model is easily adapted to the case in which  $r$ and $b$ are positive real numbers that are not necessarily integers. Often in descriptions of the model in the literature, $c$ balls are added rather than one, but adding only one ball gives simpler hypergeometric formulas, and the case of adding $c$ balls at each step with an initial inventory of $r$ red and $b$ black balls is equivalent to adding one ball at each step with an initial $r/c$ red balls and $b/c$ black balls.

%
  
An easy induction argument shows that the probability that the
first $m+n$ balls are any particular
 sequence of $m$ red and $n$ black balls is
\begin{equation*}
\frac{(r)_m (b)_n}
 {\hbox to 0pt{\hss\hspace{18pt}$(r+b)_{m+n}$\hss}},
\end{equation*}
where $(u)_n=u(u+1)\cdots(u+n-1)$.

By ball $i$, we shall mean the ball chosen at the $i$th step. Then a consequence of the previous formula is that  for any disjoint sets of integers $R=\{i_1, i_2, \dots, i_m\}$ and $B=\{j_1, j_2, \dots, j_n\}$, the probability that the balls with numbers in  $R$ are red and the balls with numbers in $B$ are black depends only on $m$ and $n$, and is therefore 
\begin{equation}
\label{e-rb}
\frac{(r)_m (b)_n}
 {\hbox to 0pt{\hss\hspace{18pt}$(r+b)_{m+n}$\hss}}.
\end{equation}  

We note that this urn model is equivalent to a lattice path model in the plane, where a particle starts at the origin, and from the point $(i, j)$ it moves right with probability $(r+i)/(r+b+i+j)$ and up with probability 
$(b+j)/(r+b+i+j)$. Then the probability that any particular path ending at $(m,n)$ has been taken is $(r)_m(b)_n/(r+b)_{m+n}$.

Now let $E_i$ be the event that ball $i$ is red, so $\E_i$ is the event that ball $i$ is black. Then with the notation of Section \ref{s-3}, for any set $S$ of positive integers, $\axx(S)$ is the probability that all the balls with numbers in $S$ are red, and 
\[\axx(S) = \frac{(r)_{|S|} }{(r+b)_{|S|}}.\]
Similarly, $\bxx(S)$ is the probability that ball $i$ is black for all $i$ in $S$, and 
\[\bxx(S)= \frac{(b)_{|S|} }{(r+b)_{|S|}}.\]
So if we take $S$ to be a set of size $n$, Theorem \ref{t-3} gives
\[\frac{(b)_n}{(b+r)_n} = \sum_{k=0}^n (-1)^k\binom nk \frac{(r)_k}{(b+r)_k},
\]
and the same identity with $b$ and $r$ switched.

This identity is a form of the Chu-Vandermonde summation theorem.
The corresponding exponential generating function identity,
\[
\egf{\frac{(b)_n}{(b+r)_n}} = e^x \egf{(-1)^n\frac{(r)_k}{(b+r)_k}},
\]
is the well-known ${}_1F_1$ transformation
\[{}_1F_1\left(\Above{b}{b+r}\Bigm | x\right) = e^x  {}_1F_1\left(\Above{r}{b+r}\Bigm | -x\right),
\]
where the hypergeometric series is defined by
\begin{equation}
{}_{p}F_q\left(\Above{u_1, \dots, u_p}{v_1,\dots, v_q } \biggm | x \right)
  = \sum_{n=0}^\infty \frac{(u_1)_n\cdots (u_p)_n}{(v_1)_n\cdots (v_q)_n} \frac{x^n}{n!}.
\end{equation}
More generally, 
we can take $m$ urns, where the $i$th urn starts with
$r_i$ red balls and  $b_i$ black balls.
At each step, we choose a ball at random from each urn, replace it, and add another ball of the same color. 
What is the probability that if we do this $n$ times, at each step we choose at least one black ball?
We let $E_i$ be the event that at the $i$th step all the balls chosen are red. Then for any finite set $S\subseteq P$, $\axx(S)$ is
the probability that for the steps in $S$, all of the balls chosen at are red. If $|S|=n$, this probability is 
\[\frac{(r_1)_n}{(r_1+b_1)_n}
  \frac{(r_2)_n}{(r_2+b_2)_n}
  \cdots
  \frac{(r_n)_n}{(r_m+b_m)_n}.
\]
So the probability that in $n$ steps at least one black ball is chosen at each step is
\begin{equation}
\label{e-urnsum}
\bxx(S) = \sum_{k=0}^n (-1)^k\binom nk 
  \frac{(r_1)_k}{(r_1+b_1)_k}
  \frac{(r_2)_k}{(r_2+b_2)_k}
  \cdots
  \frac{(r_n)_k}{(r_m+b_m)_k},
\end{equation}
for $|S|=n$,
 which may be written as the hypergeometric series
 \begin{equation}
\label{e-hyper}
{}_{m+1}F_m\left(\Above{-n,}{ }\Above{r_1,}{ r_1+b_1,}\Above{r_2,}{ r_2+b_2,}
    \Above{\cdots, }{ \cdots,}
    \Above{r_m}{ r_m+b_m} \biggm | 1\right).
\end{equation}

As a corollary, we get that this hypergeometric series is positive, as long as the $r_i$ and $b_i$ are positive real numbers.
It is not difficult to prove this result analytically, using the integral representation
\[
\frac{(r)_k}{(r+b)_{k}} = \frac{\Gamma(r+b)}{\Gamma(r)\Gamma(b)}
  \int_0^1 x^{k+r-1}(1-x)^{b-1}\, dx;
\]
However, the combinatorial approach gives a stronger result:

\begin{theorem}
\label{t-hyper}
For each nonnegative integer $n$, the hypergeometric series 
\eqref{e-hyper} can be expressed as a quotient of polynomials in the variables $r_1,\dots, r_m, b_1,
\dots b_m$ with positive coefficients.
\end{theorem}
\begin{proof}
We may assume that the variables are positive real numbers, so that the probabilistic interpretation given above applies.  Then we can compute the probability \eqref{e-urnsum} in another way by summing the probabilities of all possible outcomes in which at least one black ball is chosen at each step, and   
\eqref{e-rb} implies that each such probability is a rational function of the desired form, and so, therefore, is their sum.
\end{proof}

In the case $m=2$, the sum described in the proof of Theorem \ref{t-hyper} is simple enough to write out explicitly. We can describe the colors of the two balls chosen at each step as $(R,B)$, $(B,R)$, or $(B,B)$, where $R$ denotes red and $B$ denotes black. The number of sequences of allowable choices in which $(R,B)$ occurs $i$ times, $(B,R)$ occurs $j$ times, and $(B,B)$ occurs $k$ times is the trinomial coefficient $(i+j+k)/(i!\,j!\,k!)$. The probability of such a sequence is, by \eqref{e-rb}, 
\begin{equation*}
\frac{(r_1)_i (b_1)_{j+k}}{(r_1+b_1)_{i+j+k}}\frac{(r_2)_j (b_2)_{i+k}}{(r_2+b_2)_{i+j+k}},
\end{equation*}

and thus we have the identity
\begin{equation}
\label{e-trans}
{}_{3}F_2\left(\Above{-n,}{ }\Above{r_1,}{ r_1+b_1,}\Above{r_2}{ r_2+b_2}
     \biggm | 1\right)= \sum_{i+j+k = n} \frac{n!}{i!\,j!\,k!}
     \frac{(r_1)_i (b_1)_{j+k}}{(r_1+b_1)_{i+j+k}}\frac{(r_2)_j (b_2)_{i+k}}{(r_2+b_2)_{i+j+k}}.
\end{equation}

A slightly different approach allows us to express the double sum on the right side of 
\eqref{e-trans} as  the single sum
\begin{equation}
\label{e-singlesum}
\sum_{i=0}^n \binom ni \frac{(r_1)_i (b_1)_{n-i}(b_2)_i}{(r_1+b_1)_n(r_2+b_2)_i}.
\end{equation}
For each possible sequence of $i$ red and $n-i$ black balls chosen from the first urn, the probability of this sequence is  $(r_1)_i(b_1)_{n-i}/(r_1+b_1)_{n}$ by \eqref{e-rb}. For the choices from the second urn to be compatible, whenever a red ball is chosen from the first urn a black ball must be chosen from the second urn, but when a black ball is chosen from the first urn, the color of the ball from the second urn is unrestricted. Thus given a sequence of $i$ red and $n-i$ black balls chosen from the first urn, the probability that the choices from the second urn are compatible with those from the first urn is 
$(b_2)_i/(r_2+b_2)_i$. Summing over all possible choices for the first urn gives \eqref{e-singlesum}. 

Expressing \eqref{e-singlesum} as a hypergeometric series, we may write the identity we have proved as the  $_3F_2$ transformation
\begin{equation*}
{}_{3}F_2\left(\Above{-n,}{ }\Above{r_1,}{ r_1+b_1,}\Above{r_2}{ r_2+b_2}
     \biggm | 1\right)= 
   \frac{(b_1)_n}{(r_1+b_1)_n}
  {}_{3}F_2\left(\Above{-n,}{ }\Above{r_1,}{1-b_1-n,\,}\Above{b_2}{ r_2+b_2}
     \biggm | 1\right),
\end{equation*}
which is equivalent to formula (3.1.1) of \cite{gasper-rahman}.

%
%
A special case of the sum in \eqref{e-urnsum} deserves mention. Suppose that $\rx_1=\rx_2=\cdots=\rx_m = \rx$ and $b_1 = b_2 = \cdots = b_m=1$. Then since
\[\frac{(\rx)_k}{(\rx+1)_k} = \frac {\rx}{\rx+k},\]
the sum is 
\[U_{m,n}(\rx) = \sum_{k=0}^n (-1)^k\binom nk  
  \left(\frac {\rx}{\rx+k}\right)^m.\]

One of the interesting properties of $U_m(\rx)$ is the generating function
\[
\binom{r+n}{n}\sum_{m=0}^\infty U_{m,n}(r)\left(\frac zr\right)^m
  =\frac{z}{\prod_{i=r}^{r+n}\left(1-z/i\right)}
  \]
which is easily verified by partial fraction expansion.

These numbers (up to a constant factor) have appeared in several places in the literature. The case $r=1$ of these numbers has been studied by Smiley \cite{smiley},
who mentions a combinatorial interpretation equivalent to ours and gives further references. The case $r=2$ was considered by Foata, Han, and Lass
\cite{fhl} in connection with a coupon-collecting problem, and the general case was considered by  Laforest~\cite{laforest} and by Flajolet et al.~\cite{fgpr,flls} in the study of quadtrees. A combinatorial  connection between the coupon-collecting and quadtree interpretations and that discussed here is not apparent.

\medskip{\footnotesize
\noindent\textsc{Acknowledgment.} The author wishes to thank Arthur Benjamin for asking a question that inspired this paper, 
Christian Krattenthaler for helpful suggestions, and 
the Institut Mittag-Leffler, where this paper was written, for their hospitality.}

 \end{document}